\documentclass[11pt]{article}
\usepackage{amsmath, amsthm, amssymb, amsfonts, mathrsfs, mathtools}
\usepackage{graphicx}
\usepackage{makeidx}
\usepackage[left=2cm,right=2cm,top=2cm,bottom=2cm]{geometry}
\usepackage{caption}
\usepackage{subcaption}
\usepackage[section]{placeins}
\usepackage{enumitem}
\usepackage{tikz}
\usepackage{stmaryrd}
\usepackage{authblk}

\usetikzlibrary{decorations.pathreplacing}
\tikzstyle{vertex}=[auto=left,circle,draw=black,fill=white, inner sep=1.5]

\usepackage{hyperref}
\hypersetup{colorlinks=true, citecolor={blue}, linkcolor=blue, filecolor=magenta, urlcolor=cyan}

\newtheorem{theorem}{Theorem}[section]

\newtheorem{lema}[theorem]{Lemma}
\newtheorem{corollary}{Corollary}[theorem]

\renewcommand{\i}{\mathbf{i}}

\newcommand{\Cay}{\operatorname{Cay}}

\title{A generalization of Ramanujan's sum over finite groups}

\author[1]{Monu Kadyan}
\author[2]{Priya}
\author[3]{Sanjay Kumar Singh}
\affil[ ]{\small{\textsuperscript{1}Yau Mathematical Sciences Center, Tsinghua University, China.}}
\affil[ ]{\small{\textsuperscript{2,3}Department of Mathematics, Indian Institute of Science Education and Research Bhopal, India.}}

\affil[ ]{ {\textsuperscript{1}monu4math@gmail.com, \textsuperscript{2}priya22@iiserb.ac.in}
\textsuperscript{3}sanjayks@iiserb.ac.in}

\date{}

\begin{document}
	\maketitle
	
	\vspace{-0.3in}
	
\begin{center}{\textbf{Abstract}}\end{center} 

Let $G$ be a finite group, and let $x \in G$. Define $[x^G] := \{ y \in G : \langle x^G \rangle = \langle y^G \rangle \}$, where $\langle x^G \rangle$ denotes the normal subgroup of $G$ generated by the conjugacy class of $x$. In this paper, we determine an explicit formula for the eigenvalues of the normal Cayley graph $\Cay(G, [x^G])$. These eigenvalues can be viewed as a generalization of classical Ramanujan's sum in the setting of finite groups. Surprisingly, the formula we derive for the eigenvalues of $\Cay(G, [x^G])$ extends the known  formula of classical Ramanujan's sum to the context of finite groups. This generalization not only enrich the theory of Ramanujan’s sum but also provide new tools in spectral graph theory, representation theory, and algebraic number theory.

\vspace*{0.3cm}

\section{Introduction}
The Ramanujan's sum, introduced by Srinivasa Ramanujan in the early 20th century, is a classical arithmetic function defined over the integers. It has found profound applications in analytic number theory, particularly in the study of spectral analysis of circulant graphs. In this work, we investigate a natural and broad generalization of the Ramanujan's sum within the framework of finite groups. This generalization not only enrich the theory of Ramanujan’s sum but also provide new tools in spectral graph theory, representation theory, and algebraic number theory.

Throughout this paper, we consider $G$ to be a finite group of order $n$. We use $\textbf{1}$ to denote the identity element of $G$. Define $x^G$ as the set of all conjugates of $x$ in $G$. That is, $x^G:= \{ gxg^{-1}: g\in G \}$.   We use $\langle x^G \rangle$ to denote the normal subgroup of $G$ generated by the elements of $x^G$.

Let $\rho: G \to \mathrm{GL}_n(\mathbb{C})$ be a representation of $G$. The \textit{character} associated with $\rho$ is the function $\chi_\rho: G \to \mathbb{C}$ defined by
$$
\chi_\rho(g) := \mathrm{Tr}(\rho(g)) \quad \text{for all } g \in G,
$$
where $\mathrm{Tr}$ denotes the trace of the matrix $\rho(g)$. Throughout this paper, we adhere to the definitions and conventions established in \cite{steinberg2009representation}.

Let $S$ be a close under the inverse subset of $G$. The \textit{Cayley graph} ${\rm Cay}(G,S)$ is an undirected graph, where $V({\rm Cay}(G,S))=G$ and $$E({\rm Cay}(G,S))=\{ (a,b)\colon a,b\in G, ba^{-1} \in S \}.$$ If $S$ can be express as an union of some conjugacy classes of $G$ then ${\rm Cay}(G,S)$ is known as \textit{normal Cayley graph}. The \textit{eigenvalues} of the graph are the eigenvalues of its $(0,1)$-adjacency matrix. A graph is called \textit{integral} if all its eigenvalues are integers.

Consider the equivalence relation $\sim$ on the finite group $G$ such that $x\sim y$ if and only if $\langle x^G \rangle = \langle y^G \rangle$. For any $x\in G$, we use $[x^G]$ to denote the equivalence class of $x$ with respect to the relation $\sim$.

In 1979, Babai used characters to determine the eigenvalues of a Cayley graph. The result is as follows.
\begin{theorem}[\cite{babai1979spectra}] \label{EigenvalueExpression}
Let $G$ be a finite group, and $S$ be a close under inverse subset of $G$. If $S$ is the union of some conjugacy classes of $G$ then the eigenvalues of the Cayley graph $\Cay(G, S)$ are $$\lambda_{\chi}= \frac{1}{\chi(\textbf{1})} \sum_{s\in S} \chi (s) \textnormal{ with the multiplicity of } \chi (\textbf{1})^2,$$ where $\chi$ ranges over all irreducible characters of $G$.
\end{theorem}

Given a finite group $G$, $\chi$ be a character of $G$, and $x\in G$. Define the generalized Ramanujan sum as
\begin{align}
C_{\chi}(x):= \frac{1}{\chi(\textbf{1})} \sum_{s \in [x^G]} \chi(s).\nonumber
\end{align}

%If $\chi$ is a character of $G$, then by Theorem \ref{EigenvalueExpression}, $C_{\chi}(x)$ is an eigenvalue of $\Cay(G, [x^G])$ for each irreducible character $\chi$ of $G$. Similarly, using Theorem $1.3$ of \cite{godsil2014rationality}, $C_{\chi}(x)$ is an integer for each irreducible character $\chi$ of $G$. See~\cite{liu2012spectral} the survey paper for more information on the eigenvalues of Cayley graph. 

Let \( \chi \) be an irreducible character of a finite group \( G \). According to Theorem~\ref{EigenvalueExpression}, the value \( C_{\chi}(x) \) serves as an eigenvalue of the Cayley graph \( \mathrm{Cay}(G, [x^G]) \). Furthermore, Theorem~1.3 in \cite{godsil2014rationality} says that \( C_{\chi}(x) \) is an integer for each irreducible character \( \chi \) of \( G \). For a comprehensive overview of the spectral properties of Cayley graphs, readers are referred to the survey by Liu and Zhou~\cite{liu2012spectral}.

If $G=\mathbb{Z}_n$ and $x=1(\mod n)$, then $$C_{\chi_\alpha}(1)= \sum_{ \substack{ 1\leq s \leq n\\ \gcd(s,n)=1} } \exp \frac{2\pi \i s \alpha}{n} = R_n(\alpha),$$ where $\chi_\alpha$ is an irreducible character of $G$ such that $\chi_\alpha(s)= \exp \frac{2\pi \i s \alpha}{n} $  and $R_n(\alpha)$ is the Ramanujan's sum. It was introduced by Ramanujan in 1918. For more information, see \cite{ramanujan1918certain}. It is known that 
\begin{align}
    R_n(\alpha)=  \mu\bigg({\frac{n}{\delta_{\alpha}}}\bigg) \frac{\varphi(n)}{\varphi(\frac{n}{\delta_{\alpha}})},\label{RamanujaFormuEq}
 \end{align}
where $\delta_{\alpha} = \gcd(n,\alpha)$, $\mu$ is classic M$\ddot{\text{o}}$bius function, and $\varphi$ is Euler's phi function. Hence, we have a formula for $C_{\chi_\alpha}(1)$. Therefore, it is reasonable to think about whether a similar formula can be found for $C_{\chi}(x)$ when $G$ is a finite group, $\chi$ is any character of $G$, and $x$ is any element of $G$. 

Let $\mathcal{M}_G(\langle x^{G} \rangle)$ be the set of all the largest proper normal subgroups of $\langle x^{G} \rangle$ that are also normal subgroups of $G$. By assuming the group-theoretic version of the Chinese Remainder Theorem holds for the set $\mathcal{M}_G(\langle x^{G} \rangle)$, we obtain the following main result that generalizes to the formula of Ramanujan's sum given in Equation~(\ref{RamanujaFormuEq}). The definitions of the functions $\varphi_G$ and $\mu_G$ are provided in Section 2 and Section 3, respectively.

\noindent \textbf{Theorem} \textit{(Theorem~\ref{MainResultProof01})} Let $G$ be a finite group, $ x \in G$, and  $\chi$ be a character of $G$. If there is no proper subset $E$ of $\mathcal{M}_G(\langle x^{G} \rangle)$ with $ \bigcap\limits_{M \in E} M = \bigcap\limits_{M \in \mathcal{M}_G(\langle x^{G} \rangle)} M$, then 
\begin{align}
C_{\chi}(x) &= \mu_G(  K , \langle x^{G} \rangle )  ~ \frac{|[x^G]| }{ \varphi_G( K , \langle x^{G} \rangle ) }, \nonumber
\end{align} 
where $K=\{ x\in G : \chi(x)= \chi(\textbf{1}) \} \cap \langle x^{G} \rangle$.

The next result is a corollary specifically for abelian groups.

\noindent \textbf{Corollary} \textit{(Corollary~\ref{CoroForAbelianGroup} \cite{kadyan2024formula})} Let $G$ be a finite abelian group, $ x  \in G$, and  $\chi $ be a character of $G$. Then 
\begin{align}
C_{\chi}(x) &= \mu\bigg(\frac{\mbox{order of }x}{\mbox{size of }K}\bigg)  ~ \frac{ \varphi(\mbox{order of }x)}{ \varphi \bigg(\frac{\mbox{order of }x}{\mbox{size of }K}\bigg) }, \nonumber
\end{align}
where $K=\{ x\in G : \chi(x)= 1 \} \cap \langle x \rangle$.

%%%%%%%%%%%%%%%%%%%%%%%%%%%%%%%%%%%%%%%%%%%%%%%%%%%%%%%%%%%%%%%
%%%%%%%%%%%%%%%%%%%%%%%%%%%%%%%%%%%%%%%%%%%%%%%%%%%%%%%%%%%%%%%

\section{Preliminaries}

In this section, we present some basic definitions and notations that will be used throughout the article. The following terminology and notations will be used consistently in this paper:
\begin{itemize}
    \item  We will use the notation $\mathcal{N}(G)$ to denote the collection of all normal subgroups of $G$. 
    \item For any $V \in \mathcal{N}(G)$, define $[V^G]:= V \setminus \bigcup\limits_{\substack{U \in \mathcal{N}(G)\\ U \subsetneq V }} U$. We observe that $[x^G] = [V^G]$, where $V = \langle x^G \rangle$. If $V$ is not generated by a single conjugacy class of $G$ then $[V^G] = \emptyset$.
    
    \item Let $U \in \mathcal{N}(G)$ and $V \subseteq U$. We call $V$ to be a \textit{$G$-normal subgroup} of $U$ if $V$ is a normal subgroup of $G$.

    \item Let $U \in \mathcal{N}(G)$ and $M \subseteq U$. We call $M$ to be a \textit{maximal $G$-normal subgroup} of $U$ if it is a $G$-normal subgroup of $U$ and it is the largest $G$-normal subgroup of $U$ that is properly contained in $U$.
    
    \item For any $U \in \mathcal{N}(G)$, $\mathcal{M}_G(U)$ is the collection of all maximal $G$-normal subgroups of $U$.
    \item For any $U,V \in \mathcal{N}(G)$, $\mathcal{M}_G(U,V)$ is the collection of those maximal $G$-normal subgroups of $U$ that contain $V$.

    \item Let $U \in \mathcal{N}(G)$ and $\mathcal{M}_1, \mathcal{M}_2 \subseteq \mathcal{M}_G(U)$. We call $\mathcal{M}_1$ to be a minimal subset in $\mathcal{M}_2$ if $\mathcal{M}_1 \subseteq  \mathcal{M}_2$ and it satisfies
\begin{enumerate}[label=(\roman*)]
    \item $\bigcap\limits_{M \in \mathcal{M}_1}M = \bigcap\limits_{M \in \mathcal{M}_2}M,$
    \item $\bigcap\limits_{M \in E}M \neq \bigcap\limits_{M \in \mathcal{M}_1}M \mbox{ for all } E \subsetneq \mathcal{M}_1$.
\end{enumerate}
     If $\mathcal{M}_1$ is minimal subset in $\mathcal{M}_1$, then we call $\mathcal{M}_1$ is minimal subset in itself. Note that if $\mathcal{M}_1$ is minimal subset in $\mathcal{M}_2$, then $\mathcal{M}_1$ is minimal subset in itself, and also each subset of $\mathcal{M}_1$ is minimal subset in itself.

     \item Let $U \in \mathcal{N}(G)$ and $\mathcal{M}_G(U)$ is minimal subset in itself. For any $V \in \mathcal{N}(G)$ and $V \subseteq U$, define 
$$\varphi_G(V,U):= \left\{ \begin{array}{cl}
		1 & \mbox{if }V=U  \\
	    \prod\limits_{M \in \mathcal{M}_G(U,V)} \left( \frac{|U|}{|M|} - 1 \right)    & \mbox{if }  V  \mbox{ can be express as an intersection  of }\\ & \mbox{ some maximal $G$-normal subgroups of  } U  \\
		1   &\mbox{otherwise} 
	\end{array}\right. .$$ 
\end{itemize}

\begin{lema}\label{NewLemaForMainThm}
         Let $G$ be a finite group, $N$ be a normal subgroup of $G$, and $\mathcal{M}_G(N)$ is minimal subset in itself. If $K$ can be express as an intersection of some maximal $G$-normal subgroups of $N$ then $$\varphi_G(K,N) = \sum_{E \subseteq \mathcal{M}_G(N,K)} (-1)^{|\mathcal{M}_G(N,K)|-|E|} \prod_{M \in E} \frac{|N|}{|M|}.$$ 
\end{lema}
\begin{proof} We have
    \begin{align*}
     \varphi_G(K,N) &= \prod_{M \in \mathcal{M}_G(N,K)} \left( \frac{|N|}{|M|} - 1 \right) \\ 
     &= \sum_{E \subseteq \mathcal{M}_G(N,K)} (-1)^{|\mathcal{M}_G(N,K)|-|E|} \prod_{M \in E} \frac{|N|}{|M|}. 
    \end{align*}
   
\end{proof}

%%%%%%%%%%%%%%%%%%%%%%%%%%%%%%%%%%%%%%%%%%%%%%%%%%%%%%%%%%%
%%%%%%%%%%%%%%%%%%%%%%%%%%%%%%%%%%%%%%%%%%%%%%%%%%%%%%%%%%%%
In general, the Chinese Remainder Theorem does not hold for normal subgroups. For example, let $G=\mathbb{Z}_2 \times \mathbb{Z}_2$ and let $M_1=\langle (1,0) \rangle, M_2=\langle (0,1) \rangle$, and $M_3=\langle (1,1) \rangle$ be three maximal normal subgroups of $G$. Then $\frac{|G|}{ |M_1 \cap M_2 \cap M_3|} \neq \frac{|G|}{ |M_1|} \frac{|G|}{ |M_2|} \frac{|G|}{ |M_3|}$. However, it can be hold if we add the condition of minimality on a collection of maximal $G$-normal subgroups of $G$. The next result is an analogue to the Chinese Remainder Theorem in terms of normal subgroups.

\begin{lema}\label{LemmaMaximEquality00}
    Let $G$ be a finite group, $N$ be a normal subgroup of $G$, and $\mathcal{M}_1, \mathcal{M}_2$ be disjoint subsets of  $\mathcal{M}_G(N)$.
    \begin{enumerate}[label=(\roman*)]
        \item $\mathcal{M}_1$ is minimal subset in itself if and only if 
        \begin{align}
        \frac{|N|}{|\bigcap\limits_{M \in \mathcal{M}_1} M|} = \prod\limits_{M \in \mathcal{M}_1} \frac{|N|}{|M|}. \label{CRTeqLemma}
        \end{align}
        \item If $\mathcal{M}_G(N)$ is minimal subset in itself then $$\frac{|N|}{|\bigcap\limits_{M \in E_1 \cup E_2 } M|} = \frac{|N|}{|\bigcap\limits_{M \in E_1} M|} ~ \frac{|N|}{|\bigcap\limits_{M \in E_2} M|}$$  for all $ E_1 \subseteq \mathcal{M}_1$  and $ E_2 \subseteq \mathcal{M}_2$.
    \end{enumerate}
\end{lema}
\begin{proof}
\begin{enumerate}[label=(\roman*)]
    \item Assume that $\mathcal{M}_1$ is minimal subset in itself. Let size of $\mathcal{M}_1$ is $k$. Apply induction on $k$. If $k=1$, then statement is true. Assume that the result holds for $k=t-1$. Let $k=t$. Let $K$ be a fix element of $\mathcal{M}_1$ and $\mathcal{M}_2 := \mathcal{M}_1 \setminus \{K\}$.  Using the fact that $\mathcal{M}_1$ is minimal subset in itself, the product of both $\bigcap\limits_{M \in \mathcal{M}_2}M$ and $K$ is equal to $N$ itself. The Second Isomorphism Theorem of groups implies that 
        \begin{align}\label{IntersecEqua1}
        \frac{|N|}{|\bigcap\limits_{M \in \mathcal{M}_2} M|} = \frac{|K|}{|\bigcap\limits_{M \in \mathcal{M}_1} M|}.
        \end{align}
        By applying induction hypothesis on $\mathcal{M}_2$, Equation~(\ref{IntersecEqua1}) implies 
        \begin{align}\label{IntersecEqua2}
        \prod\limits_{M \in \mathcal{M}_2} \frac{|N|}{|M|}  = \frac{|K|}{|\bigcap\limits_{M \in \mathcal{M}_1} M|}.
        \end{align}
        Now the result follows by multiplying $|N|$ in both side in Equation~(\ref{IntersecEqua2}). Conversly, assume that Equation (\ref{CRTeqLemma}) holds. Let $\mathcal{M}_1 = \{M_1,M_2,\ldots M_k \}$. We see $$\phi: \frac{N}{\bigcap\limits_{i=1}^k M_i} \to \prod\limits_{i=1}^{k} \frac{N}{M_i}$$ with $\phi(x\bigcap\limits_{i=1}^k M_i) = (xM_1,\ldots , xM_k)$ is a group isomorphism. This isomorphism implies that the set $\mathcal{M}_1$ is minimal subset in itself.
    \item The proof follows from part $(i)$.
    \end{enumerate}
\end{proof}

%%%%%%%%%%%%%%%%%%%%%%%%%%%%%%%%%%%%%%%%%%%%%%%%%%%%%%%%%%%%%%%%%%%%%%%%%%%%%%%%%%%%%%%%%%%%%%%%%%%%%%%%%%%%%%%%%%%%%%%%%%%%%%%%%%%%%%%%%%%%%%%%%%%%%%%%%%%%%%%%%%%%%%%%%%%%%%%%%%%%%%%%%%%%%%%%%%%%%%%%%%%%%%%%%%%%%%%%%%%%%%
\begin{lema}\label{NewLemEqua22}
    Let $G$ be a finite group, and let $N,K$ are normal subgroups of $G$ with $K \subsetneq N$. If $\bigcap\limits_{M \in \mathcal{M}_G(N)} M \subseteq K$, then $K$ can be expressed as an intersection of some maximal $G$-normal subgroups of $N$. Moreover, $$K=\bigcap_{M \in \mathcal{M}_G(N,K)} M$$. 
\end{lema}
\begin{proof}
Assume that $\mathcal{M}_1$ is minimal subset in $\mathcal{M}_G(N,K)$. If $$\bigcap\limits_{M \in \mathcal{M}_1} M = \bigcap\limits_{M \in \mathcal{M}_G(N)} M$$ then $K=\bigcap\limits_{M \in \mathcal{M}_1} M$. Assume that $$\bigcap\limits_{M \in \mathcal{M}_1} M \neq \bigcap\limits_{M \in \mathcal{M}_G(N)} M.$$ Let $\mathcal{M}_2$ be a subset of $\mathcal{M}_G(N) \setminus \mathcal{M}_G(N,K)$ such that $\mathcal{M}_1 \cup \mathcal{M}_2$ is a minimal subset in $\mathcal{M}_G(N)$. We prove the result by applying induction to the size of $\mathcal{M}_2$. If the size of $\mathcal{M}_2$ is $1$, then $K=\bigcap\limits_{M \in \mathcal{M}_1} M$ because $\bigcap\limits_{M \in \mathcal{M}_1 \cup \mathcal{M}_2} M$ is maximal $G$-normal subgroup of $\bigcap\limits_{M \in \mathcal{M}_1} M$. Assume that the size of $\mathcal{M}_2$ is $t$ and $\mathcal{M}_2= \{ M_1,M_2, \ldots , M_t\}$. Define $K_1:=K \cap M_t$, $\mathcal{M}_1':=\mathcal{M}_1 \cup \{M_t\}$, and $\mathcal{M}_2':=\{ M_1,M_2,\ldots , M_{t-1} \}$. Apply induction hypothesis on $K_1$, $\mathcal{M}_1'$, and  $\mathcal{M}_2'$. We get 
$$K \cap M_t = \bigcap\limits_{M \in \mathcal{M}_1'} M= \bigcap\limits_{M \in \mathcal{M}_1} M\cap M_t.$$ The proof follows from the fact that $\bigcap\limits_{M \in \mathcal{M}_1} M\cap M_t$ is maximal $G$-normal subgroup of $\bigcap\limits_{M \in \mathcal{M}_1} M$ and $$\bigcap\limits_{M \in \mathcal{M}_1} M\cap M_t   \subseteq K \subseteq \bigcap\limits_{M \in \mathcal{M}_1} M.$$
\end{proof}

%%%%%%%%%%%%%%%%%%%%%%%%%%%%%%%%%%%%%%%%%%%%%%%%%%%%%%%%%%%%%%%%%%%%%%%%%%%%%%%%%%%%%%%%%%%%%%%%%%%%%%%%%%%%%%%%%%%%%%%%%%%%%%%%%%%%%%%%%%%%%%%%%%%%%%%%%%%%%%%%%%%%%%
%%%%%%%%%%%%%%%%%%%%%%%%%%%%%%%%%%%%%%%%%%%%%%%%%%%%%%%%%%%%%%%%%%%%%%%%%%%%%%%%%%%
%%%%%%%%%%%%%%%%%%%%%%%%%%%%%%%%%%%%%%%%%%%%%%%%%%%%%%%%%%%%%%%%%%%%%%%%%%%%%%%%%%%
%%%%%%%%%%%%%%%%%%%%%%%%%%%%%%%%%%%%%%%%%%%%%%%%%%%%%%%%%%%%%%%%%%%%%%%%%%%%%%%%%%%

\section{M$\ddot{\text{o}}$bius Inversion Formula}

The classical M$\ddot{\text{o}}$bius inversion formula is a relation between two arithmetic functions, both functions are defined from the other by taking sums over divisors. It was introduced by August Ferdinand M$\ddot{\text{o}}$bius in 1832. Its generalization was invented by Gian-Carlo Rota~\cite{rota1964foundations} in 1964. This generalization of M$\ddot{\text{o}}$bius inversion applies to functions that are defined over partially ordered sets. In 1935, Weisner \cite{weisner1935abstract} and Hall \cite{hall1936eulerian} discussed a M$\ddot{\text{o}}$bius function defined on subgroups of a group. 

In this section, we discuss a M$\ddot{\text{o}}$bius function which define on normal subgroups of $G$. Define $\mu_G : \mathcal{N}(G) \times \mathcal{N}(G) \to \mathbb{Z} $ such that 
$$\mu_G(V,U):= \left\{ \begin{array}{cl}
		1 & \mbox{if }V=U  \\
	    e-o    & \mbox{if }  V  \mbox{ can be express as an intersection  of some}\\ & \mbox{ maximal $G$-normal subgroups of  } U  \\
		0   &\mbox{otherwise} 
	\end{array}\right. $$ where $e$ (resp. $o$) is the number ways to express $V$ as the intersection of evenly (resp. oddly) many maximal $G$-normal subgroups of $U$. 

If $\mathcal{M}_G(U)$ is minimal subset in itself then 
$$\mu_G(V,U):= \left\{ \begin{array}{cl}
		1 & \mbox{if }V=U  \\
	    (-1)^r    & \mbox{if }  V  \mbox{ can be express as an intersection  of some}\\ & \mbox{ maximal $G$-normal subgroups of  } U  \\
		0   &\mbox{otherwise} 
	\end{array}\right. $$ where $M_1 , M_2 , \ldots, M_r$ are maximal $G$-normal subgroups of $U$ and $V=M_1 \cap \ldots \cap M_r$. 

    The next result give a relation between $\mu_G$ and the classical M$\ddot{\text{o}}$bius function on a finite group $G$.

\begin{lema}\label{MobiusOnAbelian}
    Let $G$ be an abelian group and $x,y \in G$. Then 
    $$\mu_G(\langle x \rangle, \langle y \rangle)=\mu\bigg(\frac{\mbox{order of } y}{\mbox{order of } x} \bigg).$$
\end{lema}
\begin{proof} We have
    $$\mu(x,y):= \left\{ \begin{array}{cl}
		1 & \mbox{if }\langle x \rangle = \langle y \rangle  \\
		(-1)^r & \mbox{if }  \langle x \rangle  = \langle p_1p_2\ldots p_r y \rangle   \\
		0   &\mbox{otherwise,} 
	\end{array}\right. $$ where $p_1, p_2, \ldots , p_r$ are some distinct prime divisors of the order of $y$. Now the result holds.
\end{proof}
    
    In general, if  $V$ can be express as an intersection  of some maximal $G$-normal subgroups of $U$, then we have 
\begin{align}
    \mu_G(V,U) = e-o  &= \sum_{\substack{ E \subseteq \mathcal{M}_G(U,V)\\ V= \bigcap\limits_{M \in E}M \\ |E| \mbox{ is even }  }} 1 - \sum_{\substack{ E \subseteq \mathcal{M}_G(U,V)\\ V= \bigcap\limits_{M \in E}M \\ |E| \mbox{ is odd }  }} 1\nonumber\\
    &= \sum_{\substack{ E \subseteq \mathcal{M}_G(U,V)\\ V= \bigcap\limits_{M \in E}M \\ |E| \mbox{ is even }  }} (-1)^{|E|} + \sum_{\substack{ E \subseteq \mathcal{M}_G(U,V)\\ V= \bigcap\limits_{M \in E}M  \\ |E| \mbox{ is odd }  }} (-1)^{|E|} \nonumber\\
    &= \sum_{\substack{ E \subseteq \mathcal{M}_G(U,V)\\ V= \bigcap\limits_{M \in E}M   }} (-1)^{|E|}.\label{MobiusDefiniEquation}
\end{align}

\begin{lema}\label{muSumSubset10} Let $G$ be a finite group. Then
$$ \sum_{\substack{ V\in \mathcal{N}( G ) \\ W \subseteq  V  \subseteq U }}  \mu_G(V,U)= \left\{ \begin{array}{cl}
		1 &  \mbox{if } W = U \\
		0   &\mbox{otherwise.} 
	\end{array}\right. $$ 	
	
\end{lema}

\begin{proof}
If $W = U$ then $$ \sum_{\substack{ V\in \mathcal{N}( G ) \\  W \subseteq  V  \subseteq U }}  \mu_G(V,U)=   \mu_G(U,U) = 1.$$ Assume that $W \subsetneq U$.  We have
$$\sum_{\substack{ V\in \mathcal{N}( G ) \\  W \subseteq  V  \subseteq U }}  \mu_G(V,U)= \sum_{\substack{ V\in \mathcal{N}( G ) \\  W \subseteq  V  \subseteq U \\ \mu_G(V,U)\neq 0 }}  \mu_G(V,U).$$ Let $M_1, M_2, \ldots , M_r$ are the only maximal $G$-normal subgroups of $U$ containing $W$. For any $0 \leq j \leq r$, ${r \choose j}$ many intersection of maximal $G$-normal subgroups exist in $\mathcal{N}( G )$. Therefore, we get
$$ \sum_{\substack{ V\in \mathcal{N}( G ) \\ W \subseteq  V  \subseteq U \\ \mu_G(V,U)\neq 0 }}  \mu_G(V,U)  = \sum_{j=0}^{r} {r \choose j} (-1)^j=(1-1)^{r}=0.$$
\end{proof}

\begin{theorem}\label{MobiusInvForm}
Let $G$ be a finite group and $f,g: \mathcal{N}(G) \to \mathbb{Z}$. If 
\begin{align}
f(U)= \sum_{\substack{ V\in \mathcal{N}(G) \\ V \subseteq U }} g(V) \label{MobiusInvForEq1}
\end{align}
then
\begin{align}
  g(U) = \sum_{\substack{ V\in \mathcal{N}(G) \\ V \subseteq U }} f(V)~ \mu_G(V,U). \label{MobiusInvForEq2}
 \end{align}
\end{theorem}
\begin{proof}
Assume that Equation (\ref{MobiusInvForEq1}) holds. We have
\begin{align}
 \sum_{\substack{ V\in \mathcal{N}(G) \\ V \subseteq U }} f(V) ~ \mu_G(V,U) &=   \sum_{\substack{ V \in \mathcal{N}( G )\\ V \subseteq U }}   \sum_{\substack{ W \in \mathcal{N}( G ) \\ W \subseteq V }} g(W) ~ \mu_G(V,U) \nonumber\\
 &=  \sum_{\substack{ W \in \mathcal{N}( G )\\ W \subseteq U }}  \sum_{ \substack{ V \in \mathcal{N}( G ) \\ W  \subseteq V \subseteq U } }    g(W)  ~ \mu_G(V,U) \nonumber\\
 &=  \sum_{\substack{ W \in \mathcal{N}( G )\\ W \subseteq U }}  g(W) \sum_{  \substack{ V \in \mathcal{N}( G ) \\ W  \subseteq V \subseteq U }  }      \mu_G(V,U) \nonumber\\
 &= g(U).\nonumber
\end{align}
Here the first equality follows from Equation (\ref{MobiusInvForEq1}) and last equality follows from Part $(i)$ of Lemma~\ref{muSumSubset10}.    
\end{proof}

The converse of Theorem~\ref{MobiusInvForm} may hold.

\begin{lema}\label{GroupSizeEqEquiLema} Let $G$ be a finite group and $ x \in G$.
\begin{align*}
|[x^{G}]| = \sum_{ \substack{ V\in \mathcal{N}(G) \\ V \subseteq \langle x^{G} \rangle } }     |V | ~ \mu_G(V,\langle x^{G} \rangle).   
\end{align*} 
\end{lema}

\begin{proof}
We have 
\begin{align*}
|\langle x^{G} \rangle | = \sum_{\substack{ V \in \mathcal{N}(G) \\ V \subseteq \langle x^{G} \rangle }} |[V^{G}]|. \nonumber
\end{align*} 
Now the result follows from Theorem~\ref{MobiusInvForm}.
\end{proof}

\begin{lema}\label{MuZeroCond}
Let $G$ be a finite group, $U,V$ are normal subgroups of $G$ and $V \subsetneq U$. If $V$ cannot be expressed as an intersection of some maximal G-normal subgroups of $U$  then $\mu_G(W,U)=0$ for all $W \in \mathcal{N}(G)$ and $W \subseteq V$.
\end{lema}
\begin{proof}
Assume that $V$ cannot be expressed as an intersection of some maximal G-normal subgroups of $U$. We get $\mu_G(V,U)=0$. If $\mu_G(W,U)\neq 0$ for some $W \in \mathcal{N}(G)$ and $W \subseteq V$, then either $W=U$ or $W$ can be expressed as an intersection of some maximal $G$-normal subgroup of $U$. If $W=U$, then $V=U$. This is not possible. On the other hand, Lemma~\ref{NewLemEqua22} implies that $V$ can be expressed as an intersection of some maximal $G$-normal subgroup of $U$. This is a contradiction.
\end{proof}

%%%%%%%%%%%%%%%%%%%%%%%%%%%%%%%%%%%%%%%%%%%%%%%%%%%%%%%%%%%%%%%
%%%%%%%%%%%%%%%%%%%%%%%%%%%%%%%%%%%%%%%%%%%%%%%%%%%%%%%%%%%%%%%
%%%%%%%%%%%%%%%%%%%%%%%%%%%%%%%%%%%%%%%%%%%%%%%%%%%%%%%%%%%%%%%

\section{Main Theorem}

Let $N$ be a normal subgroup of the finite group $G$. A character $\chi$ is \textit{principal} on $N$ if $\chi(s) = \chi(\textbf{1})$ for all $s \in N$. Define $$N^{\perp} = \{ \chi \in \widehat{G} : \chi \mbox{ is principal on } N\}.$$ 
For any character $\chi$ of $G$, define $$\mbox{Ker}(\chi)= \{ x\in G : \chi(x)= \chi(\textbf{1}) \}.$$ It is known as the kernel of the character $\chi$. That is precisely the kernel of its corresponding representation. 

\noindent For any $V \in \mathcal{N}(G)$, define $$f_{\chi}(V) := \frac{1}{\chi( \textbf{1})}  \sum\limits_{s\in V } \chi(s) \mbox{ and } g_{\chi}(V) :=  \frac{1}{\chi( \textbf{1})}  \sum\limits_{s\in [ V^{G} ] } \chi(s).$$

\noindent We have
\begin{align}
f_{\chi}(V) = \left\{ \begin{array}{cl}
		|V | &  \mbox{if }  \chi \in  V^{\perp} \\
		0   &\mbox{otherwise} 
	\end{array}\right.  \nonumber
\end{align} for all $V \in \mathcal{N}(G)$.
We can also write 
\begin{align}
f_{\chi}(\langle x^{G} \rangle) = \frac{1}{\chi( \textbf{1})} \sum\limits_{s\in \langle x^{G} \rangle } \chi(s)= \sum_{\substack{ V \in \mathcal{N}(G) \\ V \subseteq \langle x^{G} \rangle } }  \frac{1}{\chi( \textbf{1})}  \sum\limits_{s\in [ V^{G} ] } \chi(s)=  \sum_{\substack{ V \in \mathcal{N}(G) \\ V \subseteq \langle x^{G} \rangle } } g_{\chi}(V). \nonumber
\end{align}

\noindent By Theorem \ref{MobiusInvForm}, we get
\begin{align}
g_{\chi}(\langle x^{G} \rangle) &= \sum_{ \substack{ V \in \mathcal{N}(G) \\ V \subseteq \langle x^{G} \rangle } } f_{\chi}(V) ~ \mu_G(V,\langle x^{G} \rangle ) \nonumber\\
&= \sum_{ \substack{ V \in \mathcal{N}(G) \\ V \subseteq \langle x^{G} \rangle \\ \chi \in V^{\perp} } }  |V|~  \mu_G( V , \langle x^{G} \rangle ) \nonumber\\ 
&= \sum_{ \substack{ V \in \mathcal{N}(G) \\ V \subseteq \langle x^{G} \rangle \\ V \subseteq \mbox{Ker}(\chi) } }  |V|~  \mu_G( V , \langle x^{G} \rangle)\nonumber\\
&= \sum_{ \substack{ V \in \mathcal{N}(G) \\ V \subseteq  \mbox{Ker}(\chi) \cap \langle x^{G} \rangle  } }  |V|~  \mu_G( V , \langle x^{G} \rangle).\nonumber
\end{align}

\noindent The last equality follows from the fact that $\chi \in  V^{\perp}$ if and only if $V \subseteq \mbox{Ker}(\chi)$. Using the fact that $g_{\chi}(\langle x^{G} \rangle) = C_{\chi}(x)$ for all $x \in G$. We get
\begin{align}
C_{\chi}( x ) &= \sum_{ \substack{ V \in \mathcal{N}(G) \\ V \subseteq  \mbox{Ker}(\chi) \cap \langle x^{G} \rangle  } }  |V|~  \mu_G( V , \langle x^{G} \rangle).\label{SubFormulaEq1}
\end{align}

\begin{lema}\label{NewLemmaqqq} Let $G$ be a finite group, $x \in G$, and  $\chi $ be a character of $G$. Then 
\begin{align}
C_{\chi}(x) =  \sum_{ \substack{ E\subseteq \mathcal{M}_G(\langle x^{G} \rangle) \\ \bigcap\limits_{M \in E}M \subseteq  K } }  |\bigcap\limits_{M \in E}M|~~  (-1)^{|E|} , \nonumber
\end{align} 
where $K=\mbox{Ker}(\chi) \cap \langle x^{G} \rangle$. 
\end{lema}
\begin{proof}
By Equation~(\ref{SubFormulaEq1}), we get
\begin{align}
C_{\chi}( x ) &= \sum_{ \substack{ V \in \mathcal{N}(G) \\ V \subseteq  \mbox{Ker}(\chi) \cap \langle x^{G} \rangle \\ \mu_G( V , \langle x^{G} \rangle) \neq 0 } }  |V|~  \mu_G( V , \langle x^{G} \rangle)\nonumber\\
&= \sum_{ \substack{ V \in \mathcal{N}(G) \\ V \subseteq  \mbox{Ker}(\chi) \cap \langle x^{G} \rangle } }  |V|~  \sum\limits_{\substack{ E\subseteq \mathcal{M}_G(\langle x^{G} \rangle) \\  V =\bigcap\limits_{M \in E}M }} (-1)^{|E|} \nonumber\\
& = \sum_{ \substack{ E\subseteq \mathcal{M}_G(\langle x^{G} \rangle) \\ \bigcap\limits_{M \in E}M \subseteq  \mbox{Ker}(\chi) \cap \langle x^{G} \rangle } }  |\bigcap\limits_{M \in E}M|~~  (-1)^{|E|}.\nonumber
\end{align} Here, the second equality holds from Equation~(\ref{MobiusDefiniEquation}).
\end{proof}
%%%%%%%%%%%%%%%%%%%%%%%%%%%%%%%%%%%%%%%%%%%%%%%%%%%%%%%%%%%%

\begin{lema}\label{NewLemmaqqq11} Let $G$ be a finite group, $x \in G$, and  $\chi $ be a character of $G$. If $ \langle x^{G} \rangle = \mbox{Ker}(\chi) \cap \langle x^{G} \rangle$ then 
\begin{align}
C_{\chi}(x) =  |[x^{G}]|. \nonumber
\end{align} 
\end{lema}
\begin{proof}
Equation~(\ref{SubFormulaEq1}) implies that
\begin{align}
C_{\chi}(x) &= \sum_{ \substack{ V \in \mathcal{N}(G) \\ V \subseteq \langle x^{G} \rangle } }  |V|~  \mu_G( V , \langle x^{G} \rangle) \nonumber\\ 
&= |[x^{G}]|.\nonumber
\end{align}
Here the second equality holds from Lemma~\ref{GroupSizeEqEquiLema}.
\end{proof}

%%%%%%%%%%%%%%%%%%%%%%%%%%%%%%%%%%%%%%%%%%%%%%%%%%%%%%%%%%%%
%%%%%%%%%%%%%%%%%%%%%%%%%%%%%%%%%%%%%%%%%%%%%%%%%%%%%%%%%%%%

\begin{theorem}\label{MainResultProof} Let $G$ be a finite group, $x \in G$, and  $\chi $ be a character of $G$. Then 
\begin{align}
C_{\chi}(x) = \left\{ \begin{array}{cl}
		 |[x^{G}]|  &  \mbox{if } K=\langle x^{G} \rangle  \\
		\sum\limits_{ \substack{ E\subseteq \mathcal{M}_G(\langle x^{G} \rangle) \\ \bigcap\limits_{M \in E}M \subseteq  K } }  |\bigcap\limits_{M \in E}M|~  (-1)^{|E|}  & \mbox{if } K = \bigcap\limits_{M \in E}M \mbox{ for some } E\subseteq \mathcal{M}_G(\langle x^{G} \rangle) \\
        0 & \mbox{otherwise}
	\end{array}\right.   \nonumber
\end{align} 
where $K=\mbox{Ker}(\chi) \cap \langle x^{G} \rangle$. 
\end{theorem}

\begin{proof} Let $K=\mbox{Ker}(\chi) \cap \langle x^{G} \rangle$ and $\mathcal{M}= \mathcal{M}_G(\langle x^{G} \rangle) \setminus \mathcal{M}_G(\langle x^{G} \rangle,~ K )$. Rewrite the Equation~(\ref{SubFormulaEq1}) 
\begin{align}
C_{\chi}(x) =\sum_{ \substack{ V \in \mathcal{N}(G) \\ V \subseteq K  } }  |V|~  \mu_G( V , \langle x^{G} \rangle). \label{SubFormulaEq11340}
\end{align}
 We have the following three cases:

\noindent \textbf{Case 1:} Assume that $K=\langle x^{G} \rangle$. The proof follows from Lemma~\ref{NewLemmaqqq11}.\\
\noindent\textbf{Case 2:} Assume that $K \neq \langle x^{G} \rangle$ and $K$ can not be expressed as an intersection of some maximal $G$-normal subgroups of $ \langle x^{G} \rangle$. Then $\mu_G(  K , \langle x^{G} \rangle ) = 0$. Lemma~\ref{MuZeroCond} implies $\mu_G( V ,\langle x^{G} \rangle ) = 0$ for all $V \in  \mathcal{N}(G)$ and $V \subseteq  K$. Thus, Equation~(\ref{SubFormulaEq11340}) implies that $C_{\chi}(x) =0$. 
 
\noindent\textbf{Case 3:} Assume that $K$ can be expressed as an intersection of some maximal $G$-normal subgroups of $ \langle x^{G} \rangle$. Now, the proof follows from Lemma~\ref{NewLemmaqqq}. 
\end{proof}
%%%%%%%%%%%%%%%%%%%%%%%%%%%%%%%%%%%%%%%%%%%%%%%%%%%%%%%%%%%%%%%%%%%%%%%%%%%%%%%%%%%%%%%%%%%%%%%%%%%%%%%%%%%%%%%%%%%%%%%%%%%%%%%%%%%%%%%%%%%%%%%%%%%%%%%%%%%%%%%%%%%%%%%%%%%%%%%%%%%%%%%%%%%%%%%%%%%%%%%%%%%%%%%%%%%%%%%%%%%%%%

\begin{theorem}\label{MainResultProof01} Let $G$ be a finite group, $ x  \in G$, and  $\chi$ be a character of $G$. If  $\mathcal{M}_{G}(\langle x^{G} \rangle)$ is minimal subset in itself, then 
\begin{align}
C_{\chi}(x) &= \mu_G(  \mbox{Ker}(\chi) \cap \langle x^{G} \rangle , \langle x^{G} \rangle )  ~ \frac{|[x^G]| }{ \varphi_G( \mbox{Ker}(\chi) \cap \langle x^{G} \rangle , \langle x^{G} \rangle ) }. \nonumber
\end{align} 
\end{theorem}

\begin{proof}
Let $K=\mbox{Ker}(\chi) \cap \langle x^{G} \rangle$, $\mathcal{M}_1 
 = \mathcal{M}_G(\langle x^{G} \rangle,~ K )$, and $\mathcal{M}_2 = \mathcal{M}_G(\langle x^{G} \rangle) \setminus \mathcal{M}_G(\langle x^{G} \rangle,~ K )$. Rewrite the Equation~(\ref{SubFormulaEq1}) 
\begin{align}
C_{\chi}(x) =\sum_{ \substack{ V \in \mathcal{N}(G) \\ V \subseteq K  } }  |V|~  \mu_G( V , \langle x^{G} \rangle). \label{SubFormulaEq1134}
\end{align}
 
 We only need to consider the case when $K$ can be expressed as an intersection of some maximal $G$-normal subgroups of $ \langle x^{G} \rangle$. The other cases follows from Theorem \ref{MainResultProof}. We have  $$K =\bigcap\limits_{M \in \mathcal{M}_1}M.$$

 \noindent Let  $V \in \mathcal{N}(G)$, $V \subseteq K$, and $V$ can be expressed as an intersection of some maximal $G$-normal subgroups of $ \langle x^{G} \rangle$, otherwise $\mu_G(  V , \langle x^{G} \rangle )$ is zero in Equation~(\ref{SubFormulaEq1134}). We see that each $V$ in the summation of Equation~(\ref{SubFormulaEq1134}) is of this form 
\begin{align}
V = K \bigcap\limits_{ M \in E } M .
\end{align}
for some $E \subseteq \mathcal{M}_2$. Also, $\mu_G(V,\langle x^{G} \rangle) = \mu_G( K , \langle x^{G} \rangle ) ~ (-1)^{|E|}$. Apply Equation~(\ref{MobiusDefiniEquation}) in Equation~(\ref{SubFormulaEq1134}), we get
\begin{align}
C_{\chi}(x) 
&= \mu_G( K , \langle x^{G} \rangle ) \sum_{\substack{ E \subseteq \mathcal{M}_2 }} | K  \bigcap\limits_{M \in E}M|  ~ (-1)^{|E|} \nonumber\\
&= \mu_G( K , \langle x^{G} \rangle ) ~~ |\langle x^{G} \rangle| ~~ \sum_{\substack{ E \subseteq \mathcal{M}_2 }} \frac{| \bigcap\limits_{M \in \mathcal{M}_1}M  \bigcap\limits_{M \in E}M|}{|\langle x^{G} \rangle|}  ~ (-1)^{|E|} \nonumber\\
&=\mu_G( K , \langle x^{G} \rangle ) ~ |\langle x^{G} \rangle| ~ \prod_{M \in \mathcal{M}_1} \frac{|M|}{|\langle x^{G} \rangle|} \left[ 1+  \sum_{\substack{ E \subseteq \mathcal{M}_2 \\ E \neq \emptyset } } \frac{| \bigcap\limits_{M \in E} M|}{|\langle x^{G} \rangle|}   ~  (-1)^{|E|} \right].
\label{SubFormulaEq1135}
\end{align}

Here the last equality follows from Part $(ii)$ of Lemma~\ref{LemmaMaximEquality00}. By Lemma~\ref{NewLemaForMainThm}, we have 
\begin{align}
    \varphi_G( K , \langle x^{G} \rangle ) = \sum_{E \subseteq \mathcal{M}_1 } (-1)^{|\mathcal{M}_1|-|E|} \prod_{M \in E} \frac{|\langle x^G \rangle|}{|M|}.
    \label{SubFormulaEq1136}
\end{align}
\noindent Equation~(\ref{SubFormulaEq1135}) and Equation~(\ref{SubFormulaEq1136}) implies that
\begin{align}
 C_{\chi}(x) ~ \varphi_G( K , \langle x^{G} \rangle ) 
&=\mu_G( K , \langle x^{G} \rangle ) ~ |\langle x^{G} \rangle| ~ \left[ \sum_{E \subseteq \mathcal{M}_1 }  (-1)^{|\mathcal{M}_1|-|E|} \prod_{M \in E} \frac{|\langle x^G \rangle|}{|M|} \prod_{M \in \mathcal{M}_1} \frac{|M|}{|\langle x^{G} \rangle|} \right] \nonumber\\
&~~~~ \left[ 1 +   \sum_{\substack{ E \subseteq \mathcal{M}_2 \\ E \neq \emptyset } }  \frac{| \bigcap\limits_{M \in E} M|}{|\langle x^{G} \rangle|}   ~  (-1)^{|E|} \right] \nonumber\\
&=\mu_G( K , \langle x^{G} \rangle ) ~|\langle x^{G} \rangle| ~  \left[ 1 +  \sum_{ \substack{E \subseteq \mathcal{M}_1  \\ E \neq \emptyset }}  (-1)^{|E|} \prod_{M \in E} \frac{|M|}{|\langle x^G \rangle|}  \right] ~\nonumber\\
&~~~~ \left[ 1 +  \sum_{\substack{ E \subseteq \mathcal{M}_2 \\ E \neq \emptyset } } \frac{| \bigcap\limits_{M \in E} M|}{|\langle x^{G} \rangle|}    ~  (-1)^{|E|} \right] \nonumber\\
&=\mu_G( K , \langle x^{G} \rangle )  ~ |\langle x^{G} \rangle| ~\left[ 1 +  \sum_{\substack{ E \subseteq \mathcal{M}_1 \\ E \neq \emptyset }} \frac{| \bigcap\limits_{M \in E} M|}{|\langle x^G \rangle|} (-1)^{|E|}  \right] ~\nonumber\\
&~~~~ \left[ 1+  \sum_{\substack{ E \subseteq \mathcal{M}_2 \\ E \neq \emptyset } } \frac{| \bigcap\limits_{M \in E} M|}{|\langle x^{G} \rangle|}    ~  (-1)^{|E|} \right] \nonumber\\
&=\mu_G( K , \langle x^{G} \rangle ) ~ |\langle x^{G} \rangle| ~ \left[  1 + \sum_{\substack{ E \subseteq \mathcal{M}_1  \cup  \mathcal{M}_2 \\ E \neq \emptyset }}  \frac{| \bigcap\limits_{M \in E} M|}{|\langle x^{G} \rangle|}  ~  (-1)^{|E|} \right] \nonumber\\
&=\mu_G( K , \langle x^{G} \rangle ) ~ \left[  |\langle x^{G} \rangle| + \sum_{\substack{ E \subseteq \mathcal{M}_1  \cup  \mathcal{M}_2 \\ E \neq \emptyset }}  | \bigcap\limits_{M \in E} M| ~  (-1)^{|E|} \right]\nonumber\\
&=\mu_G( K , \langle x^{G} \rangle ) ~ \left[  \sum_{ \substack{ V \in \mathcal{N}(G) \\ V \subseteq \langle x^{G} \rangle  } }  |V|~  \mu_G( V , \langle x^{G} \rangle)\right] \nonumber\\
&= \mu_G( K , \langle x^{G} \rangle ) ~|[x^G]|.
\label{SubFormulaEq1137}
\end{align}
Here the fourth equality holds from Part $(ii)$ of Lemma~\ref{LemmaMaximEquality00}. Now, the proof follows from Equation~(\ref{SubFormulaEq1137}). 
\end{proof}

If $G$ is an abelian group then both $\langle x^{G} \rangle$ and $\mbox{Ker}(\chi) \cap \langle x^{G} \rangle$ are cyclic subgroups of $G$. Consequently, $\mathcal{M}_{G}(\langle x^{G} \rangle)$ is a minimal subset of itself. Based on this, we now present the following result for abelian groups.

\begin{corollary}\label{CoroForAbelianGroup}\cite{kadyan2024formula} Let $G$ be a finite abelian group, $ x  \in G$, and  $\chi $ be a character of $G$. Then 
\begin{align}
C_{\chi}(x) &= \mu\bigg(\frac{\mbox{order of }x}{\mbox{size of }K}\bigg)  ~ \frac{ \varphi(\mbox{order of }x)}{ \varphi \bigg(\frac{\mbox{order of }x}{\mbox{size of }K}\bigg) }, \nonumber
\end{align}
where $K=\mbox{Ker}(\chi) \cap \langle x^{G} \rangle$.
\end{corollary}
\begin{proof}
    By Lemma~\ref{MobiusOnAbelian}, we have $$\mu_G(K, \langle x^G \rangle )= \mu\bigg(\frac{\mbox{order of }x}{\mbox{size of }K}\bigg).$$ If this value is nonzero, i.e., $\mu_G(K, \langle x^G \rangle )\neq 0 $,   then it follows that $$\varphi_G(K, \langle x^G \rangle ) =\varphi \bigg(\frac{\mbox{order of }x}{\mbox{size of }K}\bigg).$$ Moreover, from the definition, it is clear that
    $$|[x^G]| = \varphi(\mbox{order of }x).$$ Thus, the desired result follows.
\end{proof}

\section*{Acknowledgements} The second author is supported by Junior Research Fellowship from CSIR, Government of India (File No. 09/1020(15619)/2022-EMR-I).

\noindent\textbf{Data Availability} No data was gathered or used in this paper, so a ``data availability statement'' is not applicable.

\noindent\textbf{Conflict of interest} The author states that there is no conflict of interest.
%%%%%%%%%%%%%%%%%%%%%%%%%%%%%%%%%%%%%%%%%%%%%%%%%%%%%%%%%%%%%%%%%%%%%%%%%%%%%%%%%%%%%%

%\bibliographystyle{plain}
%\bibliography{SeprateBibFile}

\end{document}